\input amstex
\documentstyle{amsppt}
\topmatter
\magnification=\magstep1
\pagewidth{5.2 in}
\pageheight{6.7 in}
\abovedisplayskip=10pt
\belowdisplayskip=10pt
\NoBlackBoxes
\title
On the $q$-extension of higher-order Euler polynomials
\endtitle
\author  Taekyun Kim  \endauthor
\affil\rm{{Division of General Education-Mathematics,}\\
{ Kwangwoon University, Seoul 139-701, S. Korea}\\
{e-mail: tkkim$\@$kw.ac.kr}}
\endaffil

\abstract{ The purpose of this paper is to present a systemic study of some families of the generalized $q$-Euler numbers
and polynomials of higher-order. In particular, by using multivariate $p$-adic invariant integral on $\Bbb Z_p$, we construct the generalized
 $q$-Euler numbers and polynomials of higher-order.
}
\endabstract
\thanks 2000 Mathematics Subject Classification  11S80, 11B68 \endthanks
\thanks Key words and phrases: Euler number, $p$-adic invariant integrals, zeta function, $p$-adic fermionic integrals \endthanks
\rightheadtext{   } \leftheadtext{ Euler numbers and
polynomials }
\endtopmatter

\document

\head\bf{ 1. Introduction }\endhead
Let $p$ be a fixed odd prime and let  $\Bbb Z_p$, $\Bbb Q_p$, $\Bbb C$ and $\Bbb C_p$ denote the ring of $p$-adic rational integers,
the field of $p$-adic rational numbers, the complex number field and the completion of algebraic closure of $\Bbb Q_p $.
The $p$-adic absolute value in $\Bbb C_p$ is normalized so that
$|p|_p=\frac1p$.
For
$d$ a fixed positive odd integer with $(p,d)=1$, let
$$X=X_d=\varprojlim_N \Bbb Z/dp^N\Bbb Z , \;\;X_1=\Bbb Z_p,$$
$$X^*=\bigcup\Sb 0<a<dp\\ (a,p)=1\endSb a+dp\Bbb Z_p,$$
$$a+dp^N\Bbb Z_p=\{x\in X\mid x\equiv a\pmod{dp^N}\},$$
where $a\in \Bbb Z$ lies in $0\leq a<dp^N$, (see [3-19]).

When one talks of $q$-extension, $q$ is variously considered as an indeterminate, a complex number $q\in\Bbb C$
or $p$-adic number $q\in \Bbb C_p$.  If $q\in\Bbb C$, one normally assumes $|q|<1$. If $q\in\Bbb C_p$, one normally assumes $|1-q|_p<1$.
In this paper we use the notation
$$[x]_q=\frac{1-q^x}{1-q}, \text{ and } [x]_{-q}=\frac{1-(-q)^x}{1+q}.$$
Let $\chi$ be the Dirichlet's character with conductor $d(=odd) \in\Bbb N$. Then the generalized Euler polynomials, $E_{n,\chi}(x)$, are defined as
$$F_{\chi}(x,t)=\frac{2\sum_{l=0}^{d-1}(-1)^l\chi(l)e^{lt}}{e^{dt}+1}e^{xt}=\sum_{l=0}^{\infty}E_{n,\chi}(x)\frac{t^n}{n!}, \text{ (see [3, 6])}. \tag1$$
We note that, by substituting $x=0$ in (1), $E_{n,\chi}(0)=E_{n,\chi}$ is the familiar $n$-th Euler number defined by
$$F_{\chi}(0,t)=\frac{2\sum_{l=0}^{d-1}(-1)^l\chi(l)e^{lt}}{e^{dt}+1}=\sum_{l=0}^{\infty}E_{n,\chi}\frac{t^n}{n!}.$$
For $ f \in UD( \Bbb Z_p) $, let us start with the expression
$$ \sum_{ 0 \leq j < p^N } (-1)^j f(j) = \sum_{0 \leq j < p^N} f(j)\mu ( j + p^N \Bbb Z_p )$$
representing analogue of Riemann's sums for $f$, cf.[1-10].

The fermionic $p$-adic invariant integral of $f$ on $\Bbb Z_p$ will be
defined as the limit $( N \rightarrow \infty )$ of these sums,
which it exists. The fermionic $p$-adic invariant integral of a function
$ f \in UD(\Bbb Z_p)$ is defined in [1, 3, 5, 7, 10] as follows:
$$I(f)= \int_{\Bbb Z_p} f(x) d \mu(x)=\lim_{N\rightarrow \infty}\sum_{0 \leq j < p^N} f(j)\mu ( j + p^N \Bbb Z_p ) = \lim_{N \rightarrow \infty}
\sum_{0 \leq j < p^N } f(j) (-1)^j .\tag2$$
Thus, we have
$$I(f_1) + I(f)=2f(0), \text{ where $f_1(x)=f(x+1)$}. $$
By using integral iterative method, we also easily see that
$$I(f_n)+(-1)^{n-1}I(f)=2\sum_{l=0}^{n-1}(-1)^{n-1-l}f(l), \text{ where $f_n(x)=f(x+n)$ for  $n \in \Bbb N$}.\tag3 $$
From (3), we note that
$$\int_{X}\chi(x)e^{xt}d\mu(x)=\frac{2\sum_{l=0}^{d-1}(-1)^l e^{lt}\chi(l)}{e^{dt}+1} =\sum_{n=0}^{\infty}E_{n,\chi}\frac{t^n}{n!}.\tag4$$
By (4), we see that
$$ \int_{X} \chi(x) x^n d\mu(x)=E_{n, \chi}, \text{ and } \int_{X}\chi(y)(x+y)^n d\mu(y)=E_{n,\chi}(x), \text{ (see [6])}. \tag5$$
The $n$-th generalized Euler polynomials of order $k$, $E_{n, \chi}^{(k)}(x)$, are defined as
$$\left( \frac{2\sum_{l=0}^{d-1}(-1)^l \chi(l) e^{lt}}{e^{dt}+1} \right)^k e^{xt}=\sum_{n=0}^{\infty} E_{n,\chi}^{(k)}(x)
\frac{t^n}{n!}, \text { (see [6, 7])}. \tag6$$
In the special case $x=0$, $E_{n,\chi}^{(k)}(0)=E_{n,\chi}^{(k)}$ are called the $n$-th generalized Euler numbers of order $k$.
Now, we consider the multivariate $p$-adic invariant integral on $\Bbb Z_p$ as follows:
$$\aligned
&\int_{X}\cdots \int_X \chi(x_1)\cdots \chi(x_k)e^{(x_1+\cdots +x_k+x)t}d\mu(x_1)\cdots d\mu(x_k)\\
&=\left(\frac{2\sum_{l=0}^{d-1}(-1)^l \chi(l)e^{lt}}{e^{dt}+1} \right)^k e^{xt}=\sum_{n=0}^{\infty}E_{n,\chi}^{(k)}(x)\frac{t^n}{n!}.
\endaligned\tag7$$
By (6) and (7), we obtain the Witt's formula for the $n$-th generalized Euler polynomials of order $k$ as follows:
$$\int_X \cdots \int_X \left(\prod_{i=1}^k \chi(x_i) \right)(x_1+\cdots +x_k +x)^n d\mu(x_1) \cdots d\mu(x_k)=E_{n,\chi}^{(k)}(x).\tag8$$
In the viewpoint of the $q$-extension of (8), we will consider the $q$-extension of generalized Euler numbers and polynomials of order $k$.
The purpose of this paper is to present a systemic study of some families of the generalized $q$-Euler numbers
and polynomials of higher-order. In particular, by using multivariate $p$-adic invariant integral on $\Bbb Z_p$, we construct the generalized
 $q$-Euler numbers and polynomials of higher-order.

\vskip 20pt

\head \bf{2. On the $q$-extension of higher-order Euler numbers and polynomials }
\endhead
\vskip 10pt
In this section we assume that $q\in\Bbb C_p$ with $|1-q|_p<1$. For $d\in\Bbb N$ with $d\equiv 1$ $(mod  \ 2)$, let $\chi$
be the Dirichlet's character with conductor $d$. For $h\in \Bbb Z, k \in \Bbb N$,
let us consider the generalized $q$-Euler numbers and polynomials of order $k$ in the viewpoint of the $q$-extension of (8).
First, we consider the $q$-extension of (1) as follows:
$$ \sum_{n=0}^{\infty} E_{n, \chi,q}(x) \frac{t^n}{n!}
=\int_X e^{[x+y]_qt}\chi(y) d\mu(y) =2\sum_{m=0}^{\infty}\chi(m)(-1)^m e^{[m]_qt}, \text{ (cf. [1, 4])}. \tag9$$
By (9), we have
$$\aligned
\int_X[x+y]_q^n\chi(y)d\mu(y)&=2\sum_{m=0}^{\infty}\chi(m)(-1)^m[m]_q^n\\
&=2\sum_{a=0}^{d-1}\chi(a)(-1)^a\frac{1}{(1-q)^n}\sum_{l=0}^n
\binom{n}{l}(-1)^l \frac{q^{l(a+x)}}{1+q^{ld}}.
\endaligned\tag10$$
From the multivariate $p$-adic invariant integral on $\Bbb Z_p$, we can also derive  the $q$-extension of
the generalized Euler polynomials of order $k$ as follows:
$$\aligned
&E_{n,\chi, q}^{(k)}(x)=\int_X \cdots \int_{X} \left(\prod_{i=1}^k \chi(x_i)\right)[x_1+\cdots +x_k+x]_q^n d\mu(x_1)\cdots d\mu(x_k)\\
&=\sum_{a_1, \cdots, a_k=0}^{d-1}\left(\prod_{i=1}^k\chi(a_i)\right)(-1)^{\sum_{j=1}^k a_j}\frac{2^k}{(1-q)^n}\sum_{l=0}^n \frac{\binom{n}{l}(-1)^lq^{l(x+\sum_{j=1}^k a_j)}}{(1+q^{dl})^k}\\
&=2^k \sum_{a_1, \cdots, a_k=0}^{d-1}\left(\prod_{i=1}^{k}\chi(a_i) \right)(-1)^{\sum_{j=1}^ka_j}\sum_{m=0}^{\infty}
\binom{m+k-1}{m}(-1)^m[x+\sum_{j=1}^k a_j+dm]_q^n\\
&=2^k \sum_{a_1, \cdots, a_{k-1}=0}^{d-1}\left(\prod_{i=1}^{k-1}\chi(a_i) \right)(-1)^{\sum_{j=1}^{k-1}a_j}\sum_{m=0}^{\infty}    
 \binom{m+k-1}{m}(-1)^m\chi(m)[x+\sum_{j=1}^{k-1} a_j+m]_q^n  .  
\endaligned\tag11$$
Let $F_{q,\chi}^{(k)}(t,x)=\sum_{n=0}^{\infty}E_{n,\chi, q}^{(k)}(x)\frac{t^n}{n!}.$ Then we have
$$\aligned
&F_{q,\chi}^{(k)}(t,x)=\sum_{n=0}^{\infty}E_{n,\chi,q}^{(k)}(x)\frac{t^n}{n!}\\
&=2^k\sum_{a_1, \cdots, a_k=0}^{d-1}\left(\prod_{i=1}^k \chi(a_i)\right)(-1)^{\sum_{j=1}^k a_j} 
\sum_{m=0}^{\infty}\binom{m+k-1}{m}(-1)^m e^{t[x+\sum_{j=1}^k a_j +dm]_q} . \endaligned\tag12   $$
From (12), we obtain the following theorem.
\proclaim{ Theorem 1} 
For $k\in \Bbb N, n\geq 0$, we have
$$\aligned
&E_{n,\chi, q}^{(k)}=\frac{2^k}{(1-q)^n}\sum_{a_1,\cdots, a_k=0}^{d-1}\left(\prod_{i=1}^k\chi(a_i)\right)(-1)^{\sum_{j=1}^k a_j}
\sum_{l=0}^n\frac{ \binom{n}{l}(-1)^lq^{l(x+\sum_{j=1}^ka_j)}}{(1+q^{ld})^k}\\
&=2^k\sum_{a_1, \cdots, a_k=0}^{d-1}\left(\prod_{i=1}^k \chi(a_i) 
\right)(-1)^{\sum_{j=1}^k a_j}\sum_{m=0}^{\infty}\binom{m+k-1}{m}(-1)^m[x+\sum_{j=1}^ka_j+md]_q^n.
\endaligned$$
\endproclaim

For $h\in\Bbb Z$, $k\in \Bbb N$, let us consider the extension of $E_{n,\chi,q}^{(k)}(x)$ as follows:
$$\aligned
 &E_{n,\chi,q}^{(h,k)}(x)=\int_X\cdots \int_X q^{\sum_{j=1}^k(h-j)x_j}
 \left(\prod_{j=1}^k \chi(x_j)\right)[x+\sum_{j=1}^kx_j]_q^n d\mu(x_1)\cdots d\mu(x_k)\\
 &=\sum_{a_1,\cdots,a_k=0}^{d-1}\left(\prod_{i=1}^k \chi(a_i)\right)(-1)^{\sum_{j=1}^ka_j}
 q^{\sum_{j=1}^ka_j(h-j)}\int_X\cdots \int_X q^{d\sum_{j=1}^k(h-j)x_j}
 \\
 &[x+\sum_{j=1}^k(dx_j+a_j)]_q^nd\mu(x_1)\cdots d\mu(x_k)\\
 &=\sum_{a_1, \cdots, a_k=0}^{d-1}\left(\prod_{i=1}^k\chi(a_i)\right)(-1)^{\sum_{j=1}^k a_j}q^{\sum_{j=1}^k a_j(h-j)}
 \frac{2^k}{(1-q)^n}\sum_{l=0}^n\frac{\binom{n}{l}(-1)^l q^{l(x+\sum_{j=1}^k a_j)}}{(-q^{d(h-k+l)}:q^d)_k},
\endaligned\tag13$$
where $( a:q)_k= (1-a) (1-aq)\cdots(1-aq^{k-1})$, (see [1, 4]).

It is well known that the Gaussian binomial coefficient is defined as
$${\binom{n}{k}}_q =\frac{ [n]_q \cdot [n-1]_q \cdots [n-k+1]_q}{ [k]_q \cdot [k-1]_q \cdots[2]_q \cdot[1]_q}, \text { (see [1, 4])}. \tag14$$
By (13) and (14), we easily see that
$$\aligned                                                                                                                              
 &E_{n,\chi,q}^{(h,k)}(x)=\frac{2^k}{(1-q)^n}\sum_{a_1, \cdots, a_k=0}^{d-1}\left(\prod_{i=1}^k \chi(a_i)\right)(-1)^{\sum_{i=1}^k a_i} 
 q^{\sum_{j=1}^k(h-j)a_j} \sum_{l=0}^n \binom{n}{l}(-1)^l \\
 &q^{l(x+\sum_{j=1}^k a_j)}\sum_{m=0}^{\infty}{\binom{m+k-1}{m}}_{q^d}(-1)^m q^{d(h-k)m}q^{dlm}\\
 &=2^k [d]_q^n \sum_{m=0}^{\infty}{\binom{m+k-1}{m}}_{q ^d}(-1)^m q^{d(h-k)m}\sum_{a_1, \cdots, a_k=0}^{d-1}
 \left(\prod_{i=1}^k\chi(a_i)\right)(-1)^{ \sum_{j=1}^k a_j}\\
 &q^{\sum_{j=1}^k (h-j)a_j}[m+\frac{x+\sum_{j=1}^k a_j}{d}]_{q^d}^n.
\endaligned\tag15$$  
Let $F_{\chi, q}^{(h,k)}(t,x)=\sum_{n=0}^{\infty}E_{n,\chi,q}^{(h,k)}(x)\frac{t^n}{n!}.$   From(15), we note that
 $$\aligned
 &F_{\chi,q}^{(h,k)}(t,x)=2^k\sum_{m=0}^{\infty}{\binom{m+k-1}{m}}_q(-1)^mq^{d(h-k)m}\sum_{a_1,\cdots, a_k=0}^{d-1}
 \left(\prod_{i=1}^k \chi(a_i) \right)(-1)^{\sum_{j=1}^ka_j}\\
 &q^{\sum_{j=1}^k(h-j)a_j}e^{t[md+x+\sum_{j=1}^k a_j]_q }.
 \endaligned\tag16$$     
 By (16), we obtain the following theorem.                                                                                          

\proclaim{ Theorem 2}
   For $h \in \Bbb Z$, $k\in\Bbb N,$ we have
$$\aligned
&E_{n,\chi,q}^{(h,k)}(x)=2^k[d]_q^n\sum_{m=0}^{\infty} {\binom{m+k-1}{m}}_q (-1)^mq^{d(h-k)m}\sum_{a_1,\cdots, a_k=0}^{d-1}
\left(\prod_{i=1}^k \chi(a_i)\right)\\
&(-1)^{\sum_{j=1}^k a_j}q^{\sum_{j=1}^k (h-j)a-j}[m+\frac{x+a_1 +a_2+\cdots+a_k}{d}]_{q^d}^n\\
&=\sum_{a_1, \cdots, a_k=0}^{d-1}\left(\prod_{i=1}^k\chi(a_i)\right)(-1)^{\sum_{j=1}^k a_j}q^{\sum_{j=1}^k (h-j)a_j}\frac{2^k}{(1-q)^n}
\sum_{l=0}^n \frac{\binom{n}{l}(-1)^l q^{l(x+\sum_{j=1}^k a_j)}}{(-q^{d(h-k+l)}:q^d)_k}.
\endaligned$$
 \endproclaim

For $h=k$, we have
$$\aligned
&E_{n,\chi, q}^{(k,k)}(x)\\
&=\frac{2^k}{(1-q)^n}\sum_{a_1, \cdots, a_k=0}^{d-1}\left(\prod_{i=1}^k\chi(a_i)\right)(-1)^{\sum_{j=1}^ka_j}
q^{\sum_{j=1}^k(h-j)a_j}
\sum_{l=0}^n\frac{\binom{n}{l}(-1)^l q^{l(\sum_{j=1}^k a_j +x)}}{(-q^{ld}:q^d)_k }\\
&=2^k [d]_q^n \sum_{m=0}^{\infty}{\binom{m+k-1}{m}}_q(-1)^m \sum_{a_1, \cdots, a_k=0}^{d-1}
\left(\prod_{i=1}^k\chi(a_i)\right)(-1)^{\sum_{j=1}^k a_j}q^{\sum_{j=1}^k(k-j)a_j}\\
&\cdot [m+\frac{x+a_1+a_2 +\cdots +a_k}{d}]_{q^d}^n .
\endaligned\tag17$$
It is not difficult to show that
$$\aligned
&\int_X \cdots \int_X \left(\prod_{j=1}^k \chi(x_j)\right)q^{\sum_{j=1}^k(m-j)x_j+mx}
d\mu(x_1)\cdots d\mu(x_k)=\sum_{a_1, \cdots, a_k=0}^{d-1}\left(\prod_{j=1}^k \chi(a_j) \right)\\
&q^{mx+\sum_{j=1}^k (m-j)a_j}
(-1)^{\sum_{j=1}^k a_j}\int_{\Bbb Z_p}\cdots \int_{\Bbb Z_p}q^{d\sum_{j=1}^k(m-j)x_j}d\mu(x_1)\cdots d\mu(x_k)\\
&=\frac{2^k q^{mx}\sum_{a_1, \cdots, a_k=0}^{d-1}\left(\prod_{j=1}^k \chi(a_j) \right)q^{\sum_{j=1}^k(m-j)a_j}(-1)^{\sum_{j=1}^k a_j}}
{(-q^{d(m-k)}:q^d)_k}
\endaligned\tag18$$
From (18), we can derive the following  equation (19).
$$\aligned
&\frac{2^k q^{mx}\sum_{a_1, \cdots, a_k=0}^{d-1}\left(\prod_{j=1}^k \chi(a_j) \right)q^{\sum_{j=1}^k(m-j)a_j}(-1)^{\sum_{j=1}^k a_j}}   
{(-q^{d(m-k)}:q^d)_k} \\
&=\int_X \cdots \int_X\left( [x+x_1+\cdots+x_k]_q(q-1)+1\right)^mq^{-\sum_{j=1}^k jx_j}\left(\prod_{j=1}^k\chi(x_j)\right)
d\mu(x_1)\cdots d\mu(x_k)\\
&=\sum_{l=0}^m\binom{m}{l}(q-1)^l \int_X \cdots \int_{X}\left(\prod_{j=1}^k\chi(x_j)\right)
[x+x_1+\cdots+x_k]_q^lq^{-\sum_{j=1}^kjx_j}d\mu(x_1)\cdots d\mu(x_k)\\
&=\sum_{l=0}^m \binom{m}{l}(q-1)^l E_{l,\chi,q}^{(0,k)}(x).
\endaligned\tag19$$
By (19), we obtain the following theorem.
\proclaim{Theorem 3}
For $d, k \in \Bbb N$ with $d\equiv 1$ $(mod  \    2)$, we have
$$\frac{2^k q^{mx}\sum_{a_1, \cdots, a_k=0}^{d-1}\left(\prod_{j=1}^k \chi(a_j) \right)q^{\sum_{j=1}^k(m-j)a_j}(-1)^{\sum_{j=1}^k a_j}} 
{(-q^{d(m-k)}:q^d)_k} =\sum_{l=0}^m \binom{m}{l}(q-1)^l E_{l,\chi,q}^{(0,k)}(x).$$
\endproclaim
From the definition of $p$-adic invariant integral on $\Bbb Z_p$, we note that
 $$\aligned
 &q^{d(h-1)}\int_X\cdots \int_X [x+d+x_1+\cdots +x_k]_q^n q^{\sum_{j=1}^k (k-j)x_j}\left(\prod_{j=1}^k \chi(x_j)\right)d\mu(x_1)\cdots d\mu(x_k)\\
 &=-\int_X \cdots \int_X [x+x_1+\cdots+x_k]_q^nq^{\sum_{j=1}^k(k-j)x_j}\left(\prod_{j=1}^k\chi(x_j)\right)d\mu(x_1)\cdots d\mu(x_k)
 +2\sum_{l=0}^{d-1}\chi(l)\\
& (-1)^l \int_X\cdots \int_X [x+\sum_{j=1}^{k-1}x_{j+1}]_q^n\left(\prod_{j=1}^{k-1}\chi(x_{j+1})\right)
 q^{\sum_{j=1}^{k-1}(h-1-j)x_{j+1}}
 d\mu(x_2)\cdots d\mu(x_k).
 \endaligned\tag20$$
By (20), we obtain the following theorem.
\proclaim{ Theorem 4}
For $h\in\Bbb Z$, $d\in\Bbb N$ with $d\equiv 1$ $(mod   \   2)$, we have
 $$q^{d(h-1)}E_{n,\chi,q}^{(h,k)}(x+d)+ E_{n,\chi,q}^{(h,k)}(x)=2\sum_{l=0}^{d-1}\chi(l)(-1)^l E_{n,q}^{(h-1,k-1)}(x). \tag21$$
Moreover, 
$$q^x E_{n,\chi,q}^{(h+1, k)}(x)=(q-1)E_{n+1, \chi,q}^{(h,k)}(x)+ E_{n,\chi,q}^{(h,k)}(x).$$
\endproclaim  
Let
$$F_{\chi, q}^{(h,1)}(t,x)=\sum_{n=0}^{\infty}E_{n,\chi,q}^{(h,1)}(x) \frac{t^n}{n!}.$$
Then we have
$$ F_{\chi, q}^{(h,1)}(t,x)=2\sum_{n=0}^{\infty}\chi(n)q^{(h-1)n}(-1)^n e^{[n+x]_qt}.\tag22$$
By (22), we see that
$$\aligned
E_{n,\chi,q}^{(h,1)}(x)&=2\sum_{m=0}^{\infty}\chi(m)q^{(h-1)m}(-1)^m[m+x]_q^n\\
&=\frac{2}{(1-q)^n}\sum_{a_1=0}^{d-1}\chi(a_1)(-1)^{a_1}\sum_{l=0}^{d-1}
\frac{\binom{n}{l}(-1)^lq^{l(x+a_1)}}{(1+q^{ld})}.
\endaligned$$


\Refs 
  \ref \no 1 \by A. Aral, V. Gupta \pages doi:10.1016/j.na.2009.07.052                                                                                         
  \paper On the Durrmeyer type modification of the $q$-Baskakov type operators \yr 2009 \vol \jour Nonlinear Analysis  \endref                     
                                                                                                                                                             
 \ref \no 2 \by M. Acikgoz, Y. Simsek \pages Art. ID 382574, P.14                                                                                                                       
 \paper On multiple interpolation functions of $N\ddot{o}rlund$-type $q$-Euler polynomials  \yr 2009 \vol 2009\jour Abst. Appl. Anal.\endref                              
                       
  \ref \no 3 \by M. Can, M. Cenkci, V. Kurt, Y. Simsek \pages 135-160                                                                                                                        
  \paper Twisted Dedekind type sums associated with Barnes' type multiple Frobenius-Euler $l$-functions,
    \yr 2009 \vol 18\jour Adv. Stud. Contemp. Math.\endref                         
                 
  \ref \no 4 \by M. Cenkci\pages 49-68                                                                                                                      
 \paper The $p$-adic  generalized twisted $(h,q)$-Euler-$l$-function and its applications \yr 2007 \vol 14\jour Adv. Stud. Contemp. Math.\endref                        
              
   \ref \no 5 \by N. K. Govil, V. Gupta \pages 97-108                                                                                                                     
  \paper Convergence of $q$-$Meyer$-$K\ddot{o}nig$-$Zeller$-$Durrmeyer$ operators  \yr 2009 \vol 19\jour Adv. Stud. Contemp. Math.\endref                          
                 
    \ref \no 6 \by T. Kim \pages 288-299                                                                                                       
   \paper $q$-Volkenborn integration  \yr 2002 \vol 9\jour Russ. J. Math. Phys.\endref                    
                     
 
   \ref \no 7 \by T. Kim \pages 151-155
\paper Symmetry identities for the twisted generalized Euler polynomials  \yr 2009 \vol 19\jour Adv. Stud. Contemp. Math.\endref

 \ref \no 8 \by T. Kim \pages 93-96 \paper Symmetry of power sum polynomials and muktivariate fermionuc $p$-adic integral on $\Bbb Z_p$                          
 \yr 2009 \vol 16 \jour Russian J. Math. Phys.                                                                                                                 
 \endref                                                                                                                                                            
                       
  \ref \no 9 \by V. Kurt \pages 2757-2764 \paper A further symmetric relation on the analogue of the Apostol-Bernoulli and the analogue of 
  the Apostol-Genocchi polynomials                           
  \yr 2009 \vol 3, no.56 \jour Appl. Math. Sciences.                                                                                                                        
  \endref                                                                                                                                                                   

 \ref \no 10 \by Y. H. Kim, K.-W. Hwang \pages 127-133                                                                                       
 \paper Symmetry of power sum and twisted Bernoulli polynomials \yr 2009 \vol 18\jour Adv. Stud. Contemp. Math.\endref                        
                     
 \ref \no 11 \by M. Cenkci  \pages 37-47 \paper The $p$-adic
generalized twisted $(h,q)$-Euler-$l$-function and its applications
\yr 2007 \vol 15 \jour Adv. Stud. Contemp. Math.\endref

 \ref \no 12 \by Z. Zhang, H. Yang \pages 191-198
\paper Some closed formulas for generalizations of Bernoulli and Euler numbers and polynomials
  \yr 2008 \vol 11\jour Proceedings of the Jangjeon Mathematical Society\endref

 \ref \no 13 \by H. Ozden, I. N. Cangul, Y. Simsek \pages 41-48
\paper Remarks on $q$-Bernoulli numbers associated with Daehee numbers
\yr 2009 \vol 18\jour Adv. Stud. Contemp. Math.\endref

 \ref \no 14 \by C.-P. Chen, L. Lin \pages 105-107
\paper An inequality for the generalized Euler constant function
  \yr 2008\vol 17\jour Adv. Stud. Contemp. Math.\endref

  \ref \no 15 \by C. S. Ryoo \pages 147-159
\paper Calculating zeros of the twisted Genocchi polynomials  \yr 2008 \vol 17
\jour Adv. Stud. Contemp. Math.\endref

 \ref \no 16 \by M. Cenkci, Y. Simsek, V. Kurt \pages 447-459
\paper Multiple two-variable $p$-adic $q$-$L$-function and its behavior at $s=0$
  \yr 2008\vol 15\jour Russ. J. Math. Phys.\endref

 \ref \no 17 \by Y. Simsek \pages 340-348 \paper On $p$-adic twisted
$q\text{-}L$-functions related to generalized twisted Bernoulli
numbers\yr 2006 \vol 13 \jour Russian J. Math. Phys.
\endref
\ref \no 18\by M. Cenkci, M. Can \pages 213-223 \paper Some results
on $q$-analogue of the Lerch zeta function \yr 2006 \vol 12\jour
Adv. Stud. Contemp. Math.
\endref

\ref \no 19\by M. Cenkci, M. Can, V. Kurt \pages 203-216 \paper
$p$-adic interpolation functions and Kummer-type congruences for
$q$-twisted and $q$-generalized twisted Euler numbers \yr 2004 \vol
9\jour Adv. Stud. Contemp. Math.
\endref \vskip 0.3cm


\endRefs
\enddocument